\documentclass[11pt]{article}
\usepackage{amssymb}
\usepackage{latexsym}

\newtheorem{thm}{Theorem}[section]
\newtheorem{prop}[thm]{Proposition}
\newtheorem{lemma}[thm]{Lemma}

\newtheorem{dfn}[thm]{Definition}

\newtheorem{rmk}[thm]{Remark}

\newcommand{\reals}{\mathbb R}
\newcommand{\complex}{\mathbb C}

\newcommand{\call}{{\cal L}}

\newcommand{\cinf}{C^{\infty}}

\newcommand{\f}{\varphi}

\def\qed{\rule{2.3mm}{2.3mm}}

\newcommand{\lag}{\langle}
\newcommand{\rg}{\rangle}

\begin{document}

\title{\bf Contact manifolds and generalized complex structures}

\author{
David Iglesias-Ponte and A\"{\i}ssa Wade\\
 {\small Department of Mathematics, The Pennsylvania State University} \\
  {\small University Park, PA 16802.} \\
  {\small e-mail: iglesias@math.psu.edu and wade@math.psu.edu}
}

\date{}
\maketitle

\begin{abstract}
We give  simple characterizations of contact $1$-forms in terms of
 Dirac structures. We  also relate normal almost contact structures
 to the theory of Dirac structures.
\end{abstract}




\footnotetext{Research partially supported by MCYT grant BFM 2003-01319.}
\section{Introduction}
 Dirac structures on  manifolds provide a unifying framework
 for the study of many geometric structures such as
  Poisson structures and closed 2-forms. They
 have applications to
  modeling of mechanical and electrical systems (see, for instance,
  \cite{BC97}).  Dirac structures were
  introduced by Courant and Weinstein (see~\cite{CW88} and \cite{C90}).
  Later, the theory of Dirac structures and Courant algebroids  was
 developed in \cite{LWX97}.

 \medskip
In \cite{Hi03}, Hitchin defined the notion of a  generalized complex
 structure on an even-dimensional manifold $M$,
 extending the setting of Dirac structures to the  complex vector bundle
 $(TM \oplus T^*M) \otimes \complex$.
  This allows to
  include other geometric structures such as  Calabi-Yau structures
 in the theory of Dirac structures. Furthermore,
 one gets a new way to look at K\"{a}hler structures
 (see  \cite{G03}).
 However, the odd-dimensional analogue of the concept
 of a  generalized complex
 structure was still missing. The aim of this Note is to fill this gap.

\medskip
 The first part of this paper concerns  
 characterizations of contact 1-forms using the notion of 
 an ${\cal E}^1(M)$-Dirac structure as introduced in~\cite{Wa00}.
 In the second part, we  define and study the  odd-dimensional analogue
 of a generalized complex structure, which
  includes the class of almost contact structures. There are 
  many distinguished subclasses of almost contact structures:
  contact metric, Sasakian, $K$-contact structures, etc. 
 We hope that the theory of Dirac structures will lead
 to  new insights on these structures.

\medskip

\section {${\cal E}^1(M)$-Dirac structures}
\subsection{Definition and examples}
In this Section, we recall the description of several
geometric structures (e.g. contact structures) in terms of Dirac
structures.

First of all,  observe that 
there is a natural bilinear form $ \lag \cdot ,\cdot
\rg$ on the vector bundle ${\cal E}^1(M)= (TM \times \reals)
\oplus (T^*M \times \reals)$ defined by:
\[
\Big \lag(X_1, f_1)+ (\alpha_1, g_1), (X_2, f_2)+ (\alpha_2, g_2)
\Big \rg =\frac{1}{2} (i_{X_2} \alpha_1 + i_{X_1} \alpha_2 +f_1 g_2 + f_2 g_1)
\]
for any $(X_j, f_j)+ (\alpha_j, g_j) \in \Gamma({\cal E}^1(M))$,
with $j=1,2$. Moreover, for any integer $k \geq 1$, one can define
$$\widetilde{d} : \Omega^k(M) \times \Omega^{k-1}(M)
 \rightarrow \Omega^{k+1}(M) \times \Omega^k(M)$$
\noindent   by the formula
$$\widetilde{d} (\alpha, \beta)= (d \alpha , \ (-1)^k \alpha + d \beta),$$
\noindent for any $\alpha \in \Omega^k(M)$, $\beta \in \Omega^{k-1}(M)$,
 where $d$ is the  exterior differentiation operator.
 When $k=0$, we define $\widetilde{d}f=(df,f)$.
 Clearly, $\widetilde{d}^2=0$.  We also  have the contraction map
  given by
$$i_{(X,f)} (\alpha , \beta)= (i_X \alpha + (-1)^{k+1} f \beta, \
  i_X \beta), $$
\noindent for  any $X \in \chi(M)$, $f \in \cinf(M)$, 
  $\alpha \in \Omega^k(M)$, $\beta \in \Omega^{k-1}(M)$.
 From these two operations, we get 
$$ \widetilde{\call}_{(X,f)}=
 i_{(X,f)} \circ \widetilde{d} + \widetilde{d} \circ i_{(X,f)}.$$
   On the space of smooth sections
of ${\cal E}^1(M)$,  we define  an operation similar to the
Courant bracket  by setting
\begin{eqnarray}\label{extended-Courant}
\kern-10pt[(X_1, f_1)+(\alpha_1, g_1),&(X_2, f_2)+(\alpha_2, g_2)
=
   ([ X_1, X_2], \ X_1\cdot f_2-X_2 \cdot f_1) \cr
  & + \widetilde{\call}_{(X_1,f_1)} (\alpha _2 ,g_2) 
- i_{(X_2 ,f_2)} \widetilde{d} (\alpha _1,g_1),
\end{eqnarray}

\noindent for any  $(X_j,f_j)+(\alpha_j,g_j) \in \Gamma({\cal E}^1(M))$ with
$j=1,2$. The skew-symmetric version of $[ \cdot , \cdot]_{{\cal E}^1(M)}$
 was introduced in \cite{Wa00}. One can notice that 
 $\widetilde{d}$ is nothing but the operator $d^{(0,1)}$ introduced 
 \cite{IM01}. Moreover,
  ${\cal E}^1(M)$ is an example of the so-called Courant-Jacobi
algebroid (see \cite{GM03}).

\begin{dfn}{\rm \cite{Wa00}}
An {\em ${\cal E}^1(M)$-Dirac structure} is a sub-bundle $L$ of
$\kern-1.5pt{\cal E}^1(M)\kern-1.5pt$  which is maximally
isotropic with respect to $ \lag \cdot,\cdot  \rg$  and
integrable, i.e., $\Gamma(L)$  is closed under the bracket $[
\cdot , \cdot ]$.
\end{dfn}

Now, we consider some examples of
${\cal E}^1(M)$-Dirac structures.

\bigskip {\bf (i) Jacobi structures}

A {\it Jacobi structure}  on a manifold $M$ is given by a pair
$(\pi, E)$ formed by a bivector field $\pi$ and a vector field $E$
such that \cite{L78}
\[
[E,\pi]_s=0, \quad [\pi, \pi]_s=2E \wedge \pi,
\]
where $[ \ , \ ]_s$ is the Schouten-Nijenhuis bracket on the space
of multi-vector fields. A manifold endowed with a Jacobi structure
is called a {\it Jacobi manifold}.  When $E$ is zero, we get a
Poisson structure.

Let $(\pi, E)$ be a pair consisting of  a bivector field  $\pi$
and a vector field $E$ on $M$. Define the bundle map $(\pi
,E)^\sharp$: $T^*M \times \reals \rightarrow TM \times \reals$  by
setting
\[
(\pi ,E)^\sharp (\alpha, g)= (\pi ^\sharp (\alpha )+ g E,
-i_{E}\alpha),
\]
where $\alpha$  is a 1-form and $g \in \cinf(M)$. The graph $L_{
(\pi, E)}$ of $(\pi ,E)^\sharp$ is an ${\cal E}^1(M)$-Dirac
structure if and only if $(\pi,E)$ is a Jacobi structure \cite{Wa00}.

\bigskip
{\bf (ii) Differential 1-forms}

Any pair $(\omega, \eta)$ formed by a 2-form $\omega$  and a
1-form $\eta$ determines a maximally isotropic sub-bundle
$L_{(\omega, \eta)}$ of ${\cal E}^1(M)$ given by
\[
(L_{(\omega, \eta)})_x =\{ (X,f)_x+(i_X \omega +f \eta, -i_X
\eta)_x: \ X \in \mathfrak X (M), \ f \in C^{\infty}(M) \}.
\]
Moreover, we have that $\Gamma (L_{(\omega ,\eta )})$ is closed
under the bracket given by (\ref{extended-Courant}) if and only if
$\omega =d \eta$. The ${\cal E}^1(M)$-Dirac structure associated 
with a 1-form $\eta$ will be denoted by $L_\eta$ (see \cite{IM02}).
\subsection{Characterization of contact structures}
In this Section, we will characterize contact structures in terms
of Dirac structures.

Let $M$ be a $(2n+1)$-dimensional smooth manifold. A 1-form $\eta$
on $M$ is {\it contact} if $\eta \wedge (d \eta)^n \ne 0$ at every
point. There arises the question of how this condition translates
into properties for $L_{\eta}$.

First, we give a characterization of Dirac
structures coming from Jacobi structures (respectively, from
differential 1-forms).
\begin{prop}\label{graph}
A sub-bundle $L$ of ${\cal E}^1(M)$ is of the form $L_{(\Lambda
,E)}$ (resp., $L_{(\omega, \eta)}$) for a pair $(\Lambda ,E)\in
\mathfrak X ^2(M)\times \mathfrak X (M)$ (resp., $(\omega
,\eta)\in \Omega ^2(M)\times \Omega ^1 (M)$) if and only if
\begin{itemize}
\item[(i)] $L$ is maximally isotropic with respect to $ \lag \cdot,\cdot
\rg$.
\item[(ii)] $L_x \cap ((T_xM \times \reals) \oplus \{0\})  =\{0\}$
(resp., $L_x \cap (\{0\}\oplus (T^\ast _xM \times \reals))
=\{0\}$) for every $x\in M$.
\end{itemize}
Moreover, $(\Lambda ,E)$ is a Jacobi structure, (resp. $\omega
=d\eta$) if and only if $\Gamma (L)$ is closed under the extended
Courant bracket (\ref{extended-Courant}).
\end{prop}

\noindent {\it Proof:} The proof of this proposition
 is straightforward (see \cite{C90} for the linear case). It is left to the reader.\hfill \qed

\medskip
Now, let $\eta$ be a contact structure on $M$. Then there exists
an isomorphism $\flat _\eta :\mathfrak X(M)\to \Omega ^1 (M)$
given by $\flat _\eta (X)= i_Xd\eta +\eta (X)\eta$ which allows us
to construct a Jacobi structure $(\pi ,E)$ given by
\[
\begin{array}{l}
\pi (\alpha ,\beta )=d\eta (\flat _\eta ^{-1}(\alpha ),
\flat _\eta ^{-1} (\beta )), \mbox{ for } \alpha ,\beta \in \Omega ^1(M),\\[5pt] E=\flat _\eta ^{-1}(\eta ),
\end{array}
\]
which satisfies that $((\pi , E)^\sharp )^{-1}(X,f)=(-i_Xd\eta -f\
\eta, \eta (X))$. Moreover, if $(\pi ,E)$ is a Jacobi structure
such that $(\pi, E)^\sharp$ is an isomorphism then it comes from a
contact structure. From these facts, we deduce that for a contact
structure $L_\eta\cong L_{(\pi ,E)}$. As a consequence of this
result and Proposition~\ref{graph}, one gets:
\begin{thm}
\label{1-to-1} There is a one-to-one correspondence between
contact 1-forms on a $(2n+1)$-dimensional manifold and ${\cal
E}^1(M)$-Dirac structures satisfying the properties
\[
\begin{array}{l}
L_x \cap ((T_xM \times \reals) \oplus \{0\})  =\{0\}, \\[8pt]
L_x \cap (\{0\}\oplus (T^\ast _xM \times \reals)) =\{0\},
\end{array}
\]
for every $x\in M$.
\end{thm}
Another characterization is the following:
\begin{thm}
An ${\cal E}^1(M)$-Dirac structure $L_\eta$ corresponds to a
contact 1-form $\eta$ if and only if
$$
L_\eta \cap ((TM \times \{0\}) \oplus (\{0\} \times \reals))
$$
is a 1-dimensional
sub-bundle of ${\cal E}^1(M)$  generated by an element of the form
$(\xi, 0)+(0, -1)$.
\end{thm}
{\it Proof:} Indeed, if $e_X=(X,0)+(0, - i_X \eta)$ then $e_X \in
L_{\eta}$ if and only if
\[
\langle (Y, g)+(i_Y d \eta+g \eta, - i_Y \eta), \ e_X \rangle=0, \
\forall \  (Y,g) \in \mathfrak X (M) \times\cinf(M),
\]
but this is equivalent to $d \eta(X,Y)=0$, for all $Y \in
\mathfrak X (M)$.

This shows $L_{\eta} \cap ((TM \times \{0\}) \oplus (\{0\}
\times \reals))$
is a 1-dimensional
sub-bundle of ${\cal E}^1(M)$ if and only if ${\rm Ker}\ d \eta$ is a
1-dimensional  sub-bundle of $TM$. If $(\xi, 0)+(0, -1)$ generates $L_{\eta} \cap (TM \times \{0\} \oplus \{0\}  \times \reals)$ then
 $$\langle (\xi, 0)+(0, -1) , \ (0,1)+(\eta,0) \rangle = \eta(\xi )-1=0.$$
Therefore,
 $${\rm Ker}\ d \eta \cap {\rm Ker}\ \eta = \{0\}.$$
We conclude that $\eta$ is a contact form.
 Moreover $\xi$ is nothing but the corresponding Reeb field, i.e., the vector field
characterized by the equations $i_\xi d\eta =0$ and $\eta (\xi
)=1$. The converse is obvious. \hfill \qed

\section{Generalized complex structures}
In this Section, we will recall the notion of generalized complex
structures.
\begin{dfn}{\rm \cite{G03}}
Let $M$ be a  smooth  even-dimensional
 manifold. A {\em generalized almost complex
structure} on $M$ is a sub-bundle $E$ of the complexification
$(TM \oplus T^*M) \otimes \complex$ such that
\begin{enumerate}
\item []{(i)} $E$ is isotropic
\item []{(ii)} $(TM  \oplus T^*M) \otimes \complex=
 E \oplus \overline{E}$, where $\overline{E}$ is the conjugate of $E$.
\end{enumerate}
\end{dfn}

The terminology  is justified by  the following result:
\begin{prop} {\rm \cite{G03}} \label{complex}
There is a one-to-one correspondence between generalized almost
complex structures and  endomorphisms ${\cal J}$ of the vector
bundle  $TM\oplus T^*M $ such that ${\cal J}^2= - id$ and ${\cal
J}$ is orthogonal with respect to $\lag \cdot , \cdot \rg$.
\end{prop}
{\it Proof:}  Suppose that $E$ is a generalized almost complex structure
on $M$. Define
\[
{\cal J} (e) = \sqrt {-1} \ e, \quad {\cal J} (\overline{e}) =-
\sqrt {-1} \ \overline{e}, \quad \mbox{for any} \  e \in
\Gamma(E).
\]
Then, ${\cal J}$ satisfies the properties ${\cal J}^2= - id$ and
${\cal J}^*= - {\cal J}$. Conversely, assume that ${\cal J}$
satisfies these two properties. Define the sub-bundle $E$ whose
fibre as the $\sqrt{-1}$-eigenspace of $\cal J$. It is not
difficult to prove that $E$ is isotropic under $\lag \cdot , \cdot
\rg$. Moreover, since $\overline{E}$ is just the
$(-\sqrt{-1})$-eigenspace of $\cal J$ we get that $(TM  \oplus
T^*M) \otimes \complex= E \oplus \overline{E}$. \hfill \qed

\medskip

We have the following definition:

\begin{dfn}
Let $M$ be an even-dimensional smooth manifold. A generalized
almost complex structure $E \subset (TM  \oplus T^*M ) \otimes
\complex$ is {\em integrable} if  it is closed under the Courant
bracket. Such a sub-bundle is called a {\em generalized complex
structure}.
\end{dfn}

The$\kern-.5pt$ notion$\kern-.5pt$ of a$\kern-.5pt$
generalized$\kern-.5pt$ complex structure$\kern-.5pt$ on
an$\kern-.5pt$ even-dimensional$\kern-.5pt$ smooth manifold was
introduced by Hitchin in \cite{Hi03}.

\section{Generalized almost contact structures}
The existence of a generalized almost complex structure on $M$
forces the dimension  of $M$ to be even (see \cite{G03}). A
natural question to ask is: what would be the odd-dimensional
analogue of a generalized almost complex structure?

To define the analogue of the concept of a generalized almost
complex structure for odd-dimensional manifolds, one should
consider the vector bundle ${\cal E}^1(M) \otimes \complex$
instead of $ (TM  \oplus T^*M) \otimes \complex$.
\begin{dfn}
Let $M$ be a real smooth manifold of dimension $d=2n+1$. A {\em
generalized almost contact structure } on $M$ is a sub-bundle $E$
of $ {\cal E}^1(M) \otimes \complex$ such that $E$ is isotropic
and
$$ {\cal E}^1(M) \otimes \complex= E \oplus \overline{E},$$
where $\overline{E}$ is the complex conjugate of $E$.
\end{dfn}
 By a proof similar
 to  that of Proposition \ref{complex}, one gets the following
result.
\begin{prop}
Let $M$ be a real smooth manifold of dimension $d=2n+1$. There is
a one-to-one correspondence between generalized almost contact
structures on  $M$ and  endomorphisms ${\cal J}$ of the vector
bundle  $ {\cal E}^1(M) $ such that ${\cal J}^2= - id$ and ${\cal
J}$ is orthogonal with respect to $\lag \cdot , \cdot \rg$.
\end{prop}

\subsection{Examples}
{\bf (i) Almost contact structures}.

Let  $M$ be  a smooth manifold of dimension $d=2n+1$. An {\em
almost contact structure}  on $M$ is a triple $(\varphi, \xi,
\eta)$, where $\varphi$ is a (1,1)-tensor field, $\xi $ is a
vector field on  $M$, and $\eta$ is a 1-form such that
$$
\eta(\xi)=1 \quad {\rm and} \quad \varphi^2(X)= - X + \eta(X) \xi, \quad 
 \forall \ X \in
 \mathfrak X(M)
$$
(see \cite{Bl02}). As a first consequence, we get that
$$
\varphi (\xi )=0, \qquad \eta \circ \varphi =0 .
$$
We now show that every almost contact structure  determines a
generalized almost contact structure. Define $J :  \Gamma(TM
\times \reals) \rightarrow \Gamma(TM \times \reals)$ by
$$
J(X, f)= (\varphi X-f \xi, \ \eta(X)), \quad \mbox{for all} \
 X \in \mathfrak X(M), f \in \cinf(M).
$$
Then $J^2 = -id$. Let $J^*$ be the dual map of $J$. Consider the
endomorphism ${\cal J}$ defined by
\[
{\cal J}(u)= J(X,f)- J^*(\alpha, g).
\]
for $u=(X,f)+(\alpha, g)\in \Gamma({\cal E}^1(M) )$. Then  ${\cal
J}$ satisfies ${\cal J}^2= -id$ and ${\cal J}^*= - {\cal J}$.

In addition, one can deduce that the generalized almost contact
structure $E$ is given by
\begin{equation}\label{bundle-almostcontact}
E= F \oplus Ann(F),
\end{equation}
where
\begin{equation}\label{bundle-almostcontact2}
F_x= \{  \ J( X,f)_x +\sqrt{-1} (X,f)_x \ | \
 (X,f) \in \Gamma(TM \times \reals) \}
\end{equation}
and $Ann(F)$ is the annihilator of $E$.

\bigskip

{\bf (ii) Almost cosymplectic structures}

An almost cosymplectic structure on a smooth manifold $M$ of dimension
$d=2n+1$ is a pair $(\omega, \eta)$ formed by a 2-form $\omega$
and a 1-form $\eta$ such that $\eta \wedge \omega^n \ne 0$ everywhere. 
The map $\flat:  \mathfrak X(M) \rightarrow
\Omega^1(M)$  defined by
$$\flat(X)= i_X\omega+ \eta(X) \eta,  \quad \forall \ X \in \mathfrak X(M).$$
is an isomorphism of $\cinf(M)$-modules. The vector field $\xi=
\flat^{-1}(\eta)$ is called the Reeb vector field of the
almost cosymplectic structure and it is characterized by
$i_{\xi}\omega=0$ and $\eta(\xi)=1$. Define $\Theta: \mathfrak
X(M) \times \cinf(M) \rightarrow \Omega^1(M) \times \cinf(M)$ by
$$\Theta(X,f)= \Big(i_X\omega+f \eta, \ -\eta(X)\Big),
 \quad \forall \ X \in \mathfrak X(M), \ \forall \ f \in \cinf(M).$$

One can check that $\Theta$ is an isomorphism of $\cinf(M)$-modules.
 Let  ${\cal J}:
 \Gamma({\cal E}^1(M)) \rightarrow \Gamma({\cal E}^1(M))$ be the endomorphism
 given by
 $${\cal J}\Big((X,f)+(\alpha, g) \Big)= -\Theta^{-1}(\alpha, g)+
 \Theta(X,f).$$
It is easy to check that ${\cal J}^2=-id$. Moreover,
 for $e_i=(X_i,f_i)+(\alpha_i, g_i) \in \Gamma({\cal E}^1(M)) $, we have
$$\lag {\cal J}e_1, \ e_2 \rg=
 \lag -\Theta^{-1}(\alpha_1, g_1)+
 \Theta(X_1,f_1), \ (X_2,f_2)+(\alpha_2, g_2)\rg =-
\lag  e_1, \ {\cal J}e_2 \rg.$$
Hence ${\cal J}^*=-{\cal J}$.

This shows that every almost cosymplectic structure determines a
generalized almost contact structure. Furthermore, the associated
bundle $E$ is given by
\begin{equation}\label{bundle-cosymplectic}
E_x= \{  \ (X,f)_x-\sqrt{-1} \Theta(X,f)_x \ | \
 (X,f) \in \Gamma(TM \times \reals) \}.
\end{equation}
\section{Integrability}
By analogy to generalized complex structures, one
can consider the integrability of a generalized almost contact
structure.
\begin{dfn}
On an odd-dimensional smooth manifold $M$, we say that a
generalized almost contact structure $E \subset {\cal E}^1(M)
\otimes \complex$ is {\em integrable} if it is closed under the
extended Courant bracket given by Eq. (\ref{extended-Courant}).
\end{dfn}
\subsection{Examples}
{\bf (i) Normal almost contact structures}

\medskip
An almost contact structure  $(\varphi, \xi, \eta)$ is {\em
normal} if
$$N_{\varphi}(X,Y)+d\eta(X,Y)\xi=0, \quad \mbox{for all} \  X, Y \in
 \mathfrak X(M), $$
where $N_{\varphi}$ is the Nijenhuis torsion of $\varphi$, i.e.,
\[
N_{\varphi}(X,Y)= [\varphi X, \varphi Y] + \varphi^2 [X, Y] -
\varphi[ \varphi X, Y]- \varphi[X, \varphi Y].
\]
Some properties of normal almost contact structures are the following 
ones (see \cite{Bl02}).
\begin{lemma}\label{Lemma-normal}
If an almost contact structure $(\varphi, \xi,\eta )$ is normal
then it follows that
\[
\begin{array}{ll}
d \eta( X, \xi)=0,& \eta [\varphi X, \xi ]=0,\\[5pt]
[\varphi X, \  \xi]= \varphi [X,  \xi] &d\eta( \varphi X,
Y) = d\eta( \varphi Y , X),
\end{array}
\]
for $X,Y \in \mathfrak X (M)$.
\end{lemma}
{\it Proof:} Applying normality condition to $Y=\xi$ we get that
\begin{eqnarray*}
0=N_{\varphi}( X,  \xi)+d \eta( X,  \xi) \xi=
  \varphi^2[ X,  \xi]- \varphi[ \varphi X,  \xi] +
  d \eta( X, \xi) \xi .
\end{eqnarray*}
Using the fact that $\eta \circ \varphi=0$, we obtain $d \eta( X,
\xi)=0,$ for any $X \in \mathfrak X(M)$. As a consequence, $\eta
[\varphi X, \xi ]=0$. On the other hand,
\begin{eqnarray*}
0&=&N_{\varphi}(\varphi X,  \xi) +d \eta(\varphi X, \xi) \xi\cr
&=& \varphi^2[\varphi X,  \xi]- \varphi[ \varphi^2X,  \xi] +
  d \eta(\varphi X, \xi) \xi \cr
&=& - [\varphi X,  \xi]+ \varphi [X,  \xi],
\end{eqnarray*}
Finally, if $X,Y\in \mathfrak X(M)$ then
\[
\eta(N_{\varphi}(\varphi X, Y) +d \eta(\varphi X,Y) \xi)=
  - \eta([ \varphi ^2 X, Y] + [\varphi X,  \varphi Y]) +d\eta (\varphi X, Y).
\]
We deduce  that $d\eta( \varphi X, Y) = d\eta( \varphi Y ,
X)$.\hfill \qed

\medskip
We have seen that every almost contact structure $(\varphi, \xi,
\eta)$ determines a generalized almost complex structure $E\subset
{\cal E}^1(M) \otimes \complex$. Furthermore, we have the
following result:
\begin{thm}
An almost contact structure $(\varphi, \xi, \eta)$ is normal if
and only if its corresponding sub-bundle $E$ given by
(\ref{bundle-almostcontact}) is integrable.
\end{thm}
{\it Proof:} Clearly, the integrability of $E$ is equivalent to
the closedness of $\Gamma(F)$ under the extended Courant bracket,
where $F$ is the sub-bundle defined by
(\ref{bundle-almostcontact2}). Suppose $[\Gamma(F), \Gamma(F)]
\subset \Gamma(F)$. Let $u_X=(X,0)$, $u_Y =(Y, 0) \in \Gamma({\cal
E}^1(M))$. Denote $e_X= Ju_X+\sqrt{-1} \ u_X$ and $e_Y=J
u_Y+\sqrt{-1} \ u_Y$. Then
$$[e_X,  e_Y] \in F \iff
[Ju_X, Ju_Y]- [u_X, u_Y]= J\Big([J u_X, u_Y]+
[u_X, Ju_Y] \Big).$$
By a simple computation, one gets
$$[Ju_X, Ju_Y]- [u_X, u_Y]=
 \Big([\varphi X, \varphi Y]-[X, Y], \ \varphi X \cdot \eta(Y)-
 \varphi Y \cdot \eta(X) \Big).$$
Moreover, the term
$J\Big([Ju_X, u_Y]+[u_X, Ju_Y] \Big)$ equals
 $$\Big(
 \varphi([ \varphi X, Y] + [X,  \varphi Y]) -
 (X\cdot \eta(Y)- Y \cdot \eta(X)) \xi, \
 \eta([ \varphi X, Y] + [X,  \varphi Y]\Big).$$
Therefore $[e_X,  e_Y] \in \Gamma(F)$ if and only if
 \begin{displaymath}\left \{
\begin{array}{ccc}
  [\varphi X, \varphi Y]-[X, Y]=\varphi([ \varphi X, Y] + [X,  \varphi Y]) -
 (X\cdot \eta(Y)- Y \cdot \eta(X)) \xi \cr
 {}\cr
 \varphi X \cdot \eta(Y)-
 \varphi Y \cdot \eta(X)=
 \eta([ \varphi X, Y] + [X,  \varphi Y])  \ \ \ \ \ \ \ \ \ \ \ \ \ \ \
\ \ \ \ \ \ \ \ \ \ \ \ \ \ \
\end{array}\right.
\end{displaymath}
Because $[X, Y]= - \varphi ^2([X,Y]) +
 \eta([X,Y])\xi$ and $\eta (\varphi X)=0$, for any $X$, $Y \in \mathfrak X(M)$,
  this implies the relations
 \begin{displaymath}\left \{
\begin{array}{ccc}
  N_{\varphi}(X, Y) +d \eta(X,Y) \xi=0 \cr
 {}\cr
 d\eta( \varphi X, Y) = d\eta( \varphi Y , X)
\end{array}\right.
\end{displaymath}
This proves that  if $E$ is integrable  then the almost contact
structure is normal. Conversely, suppose that $N_{\varphi}(X, Y)
+d \eta(X,Y) \xi=0$, for any $X, Y$ in $\mathfrak X(M)$. Using
Lemma \ref{Lemma-normal}, we also have that $d\eta( \varphi X, Y)
= d\eta( \varphi Y , X)$. Thus, we conclude that $[e_X,  e_Y] \in
\Gamma(F)$, for any $e_X= u_X+\sqrt{-1} \ Ju_X$, $e_Y=
u_Y+\sqrt{-1}\  Ju_Y$ in $\Gamma(F)$.

It remains to show that $[ e_X, \  J(0,1)+\sqrt{-1} (0,1)]$ is in
$ \Gamma(F) $, for any  section $e_X= Ju_X+\sqrt{-1}\  u_X \in
\Gamma(F) $. This condition is equivalent to the relations
\begin{displaymath}\left \{
\begin{array}{ccc}
[\varphi X, \  \xi]= \varphi [X,  \xi]  \cr
{}\cr
  \xi \cdot \eta(X) =  -\eta([X,  \  \xi])
\end{array}\right.
\end{displaymath}
The relation $ \xi \cdot \eta(X) =  -\eta([X,  \  \xi]) $ is
satisfied since $d \eta(X, \xi)=0$ by Lemma \ref{Lemma-normal}. We
conclude that $[e_X, \ J(0,1)+\sqrt{-1} (0,1)] \in F $. Therefore
$F$ is closed under that extended Courant bracket, which means that
$E$ is integrable.\hfill \qed

\bigskip
{\bf (ii) Contact structures}

\medskip
Let $(\omega ,\eta )$ be an almost cosymplectic structure and $E$ the
associated generalized almost contact structure given by
(\ref{bundle-cosymplectic}). We will prove that the integrability
condition forces $\eta$ to be a contact structure. In fact,
\begin{prop}
Let $(\omega ,\eta )$ be an almost cosymplectic structure on a manifold
$M$ and $E$ the associated generalized almost contact structure.
Then, $E$ is integrable if and only if $\omega =d\eta$. As a
consequence, $\eta$ is a contact structure on $M$.
\end{prop}
{\it Proof:} Let $e_1,e_2\in \Gamma (E)$. One can easily show that
$[e_1,e_2]\in \Gamma (E)$ if and only if $\omega =d\eta$. \hfill
\qed

\begin{rmk}
{\rm Following \cite{G03}, one can define an analogue of
generalized K\"{a}hler structure. In our setting, one could define the
notion of a {generalized Sasakian structure} as a pair $({\cal
J}_1,{\cal J}_2)$  of commuting generalized integrable generalized
almost contact structures, i.e. ${\cal J}_1\circ {\cal J}_2={\cal
J}_2\circ {\cal J}_1$, such that $G=- {\cal J}_1{\cal J}_2$
defines a positive definite metric on ${\cal E}^1(M)$. In
particular, every Sasakian structure is a generalized Sasakian
structure. We postpone  the study of this
notion and its main properties to a separate paper. }
\end{rmk}
{\bf Acknowledgments:}
 D.~Iglesias wishes to thank the Spanish Ministry of Education and Culture and Fulbright program for a MECD/Fulbright postdoctoral grant.

\end{document}